\magnification=1200
\input amstex
\documentstyle{amsppt}
\hsize 6truein
%\voffset .6truein
\vsize 8.5truein
\def\Ind{{\hbox{Ind}}}
\def\leftitem#1{\item{\hbox to\parindent{\enspace#1\hfill}}}
\def\Vir{\operatorname{Vir}}
\def\dim{\operatorname{dim}}

\def\bC{\Bbb C}

\def\Z{\Bbb Z}\def\N{\Bbb N}
\def\bQ{\Bbb Q}

\def\f{ \hbox{ \ for \ }}

   \def\b{\beta}
  \def\a{\alpha}
 \def\d{\delta} \def\g{\gamma}

\def\C{\Bbb C}
\def\al{\alpha}
\def\be{\beta}

\def\ab{{\alpha+\beta}}
\def\V{{V}}
\def\gg{{\frak g}}
\def\hh{{\frak h}}
\def\Prod{\Pi}

\def\vs{{\vskip 5pt}}
\def\tilde{\widetilde}

\leftheadtext{Billig and Zhao}

\rightheadtext{modules with exponential highest weight}

\topmatter
\title
 Weight modules over exp-polynomial Lie algebras
\endtitle
\author Yuly Billig and Kaiming Zhao
\endauthor
\affil School of Mathematics and Statistics\\
      Carleton University \\
       Ottawa, Ontario, K1S 5B6 Canada \\
              Email: billig\@math.carleton.ca\\
\\
and\\
\\
Institute of Mathematics\\
Academy of Mathematics and System Sciences\\
Chinese Academy of Sciences\\
Beijing 100080, P. R. China\\
       Email: kzhao\@mail2.math.ac.cn\\
\endaffil
\thanks
Research supported by  Natural Sciences and Engineering Research Council of Canada
(Y.B.) and  Hundred Talents Program of Chinese
Academy of Sciences (K.Z.).
\endthanks

\abstract In this paper, we generalize a result by Berman and Billig
on weight modules over Lie algebras with polynomial multiplication.
More precisely, we show that a highest weight module
with an exp-polynomial ``highest weight'' has finite dimensional weight
spaces.  We also get a class of irreducible weight modules with finite dimensional
weight spaces over generalized Virasoro algebras which do not occur over the
classical Virasoro algebra.

\endabstract

\endtopmatter
\document

\bigskip
\subhead 0. \ \ Introduction
\endsubhead
\medskip

Representations of affine Lie algebras and the Virasoro algebra
have many important applications in mathematics and physics. One
of the main ingredients of these theories is the construction of
the highest weight modules. Recently there has been substantial
activity in developing representation theories for higher rank
infinite dimensional Lie algebras, in particular toroidal Lie
algebras, generalized Virasoro algebras and quantum torus algebras (see [EM],
[GL], [L], [BS], [G], [BB], [B1], [B2], [Ma], [HWZ]).

Unlike rank one algebras (affine and Virasoro), the higher rank
infinite dimensional Lie algebras do not possess a triangular
decomposition, which makes the standard construction of the
highest weight modules inapplicable. Nonetheless, there have been
found several explicit realizations of representations for these
algebras using the vertex operator approach.

In the vertex constructions the highest weight is replaced with a
loop-like module for the subalgebra of degree zero (in general, this
subalgebra is not commutative). Let us describe in brief
these representations from the perspective of the highest weight
modules.

 Let $G$ be a $\Z \times \Z^n$-graded Lie algebra and let $G = G^- \oplus G^{(0)}
\oplus G^+$ be a decomposition of $G$ relative to $\Z$-grading.
The subalgebra $G^{(0)}$ is an infinite-dimensional Lie algebra of
rank $n$. We take some natural module $V$ for $G^{(0)}$ (usually
$V$ is either $\Z^n$-graded or finite dimensional). Parallel to
the construction of a highest weight module, we let $G^+$ act on
$V$ trivially, and introduce the induced module
$$\tilde M(V)={\roman{Ind}}^G_{G^{(0)}+G^+}V\simeq U(G^-)\otimes_\bC V.$$
If $V$ is $\Z^n$-graded then $\tilde M(V)$ inherits a $\Z \times \Z^n$-grading, and it
is $\Z$-graded when $V$ is finite-dimensional.

 The difficulty here is that $\tilde M(V)$ will have infinite-dimensional homogeneous
components (and thus will not have a character). Nonetheless the explicit vertex
operator constructions show that in some cases $\tilde M(V)$ has quotients with
finite-dimensional homogeneous components. This situation has been clarified in [BB],
where it was proved that $\tilde M(V)$ has a graded factor-module $M(V)$ with
finite-dimensional components provided that $G$ is a polynomial Lie algebra and
$V$ is a polynomial module. The polynomiality condition means that the structure
constants of the Lie algebra and of the module are given by polynomial expressions
(see Section 1 for the precise definitions). However this left out some important
examples, notably quantum torus Lie algebras. In this paper we expand the class of Lie
algebras and modules for which the theorem is applicable. We prove that $M(V)$ has
finite-dimensional homogeneous components when $G$ is an exp-polynomial Lie algebra
and $V$ is an exp-polynomial module.

 We illustrate our definitions and theorems with a sequence of 24 examples.

\smallskip
The paper is organized as follows.
In Section 1 we define exp-polynomial Lie algebras, exp-polynomial modules
and give the statement of the main results
(Theorems 1.5 and 1.7). We provide numerous examples (old and new)
of  such Lie algebras and modules as well as give examples of applications of
the main theorems.
In Section 2 we first establish the extended Vandermonde
determinant formula, and then  give the proof of Theorems 1.5 and 1.7.
In Section 3 we
show a  similar but stronger result for  generalized Virasoro
algebras to give a class of irreducible weight modules with finite
dimensional weight spaces over generalized Virasoro algebras which
do not occur over the classical Virasoro algebra.

\bigskip
\subhead 1. \ \ Exp-polynomial Lie algebras and exp-polynomial modules
\endsubhead
\medskip

Let $\bC$ be the complex number field.  We assume that all Lie algebras and vector spaces
are over $\bC$ in this paper, although $\bC$ can be replaced
by any field of characteristic $0$.

{\bf Definition 1.1.} The algebra
% $\Cal A_r$
of exp-polynomial
functions in $r$ variables, $n_1, \ldots,$
\break
$n_r$, is the algebra of
functions $f(n_1, \ldots, n_r): \Z^r \rightarrow \C$ generated as
an algebra by functions $n_j$ and $a^{n_j}$, where
$a\in\C^*=\C\setminus\{0\}, j=1,\ldots,r$.

An exp-polynomial function may be written as a finite sum
$$f(n_1, \ldots, n_r) = \sum_{k \in \Z_+^r} \sum_{a\in(\C^*)^r} c_{k,a}
n_1^{k_1} \ldots n_r^{k_r} a_1^{n_1} \ldots a_r^{n_r} ,$$ where
$c_{k,a}\in\C$, $k = (k_1,\ldots,k_r)$ with $k_j \geq 0$, and $a
= (a_1, \ldots, a_r)$.

An exp-polynomial function $f(n)$ has the property that the
function $f(n+m)$ is also exp-polynomial as a function of two
variables.

% Let us first recall the definition of  Lie algebras with polynomial multiplications defined in % [BB].

\smallskip
{\bf Definition 1.2.}  Let $G=\oplus_{\a\in\Z^n}G_\a$ be a
$\Z^n$-graded Lie algebra. Then $G$ is said to be an {\it
exp-polynomial Lie algebra} if $G$ has a homogeneous spanning set
$\{g_k(\a)\,\,\,|\,\,\,k\in K,\,\,\a\in\Z^n\}$ with
 $g_k(\a)\in G_\a$,
and there exists a family of exp-polynomial functions
$\{f_{k,r}^s(\a,\b)\,\,\,|\,\,\,k,r,s\in K\}$ in the $2n$ variables $\a_j,\b_j$
 and where for each $k,r$ the set $\{s|f_{k,r}^s(\a,\b)\ne0\}$ is finite, such that
$$[g_k(\a),g_r(\b)]=\sum_{s\in K}f_{k,r}^s(\a,\b)g_s(\a+\b), \; \; \f k,r\in K,\,
\a,\b\in\Z^n.\eqno(1.1)$$
This homogeneous spanning set $\{g_k(\a)\,\,\,|\,\,\,k\in
K,\,\,\a\in\Z^n\}$ is call {\it a distinguished spanning set}.

% Actually there is no difference between spanning set and spanning
% basis in this definition. Sometimes spanning set can help to
% easily see the functions. In later proofs, using basis is more
% convenient.
%
\smallskip
% We shall be interested in only the Lie algebras for $n\ge1$.
%
If the functions
$f_{k,r}^s(\a,\b)$ in (1.1) are in fact polynomials, we say that $G$ is a
{\it polynomial} Lie algebra (cf. [Definition 1.6, BB]).

\smallskip
{\bf Example 1.} Let $\gg$ be a finite dimensional Lie algebra
with a basis $B = \{ g_k \}_{k \in K}$. Then the toroidal Lie
algebra $G = \C[t_1^\pm, \ldots, t_n^\pm] \otimes \gg$ is a
polynomial Lie algebra with the distinguished spanning set
$g_k(\al) = t^\al g_k, \; \al \in \Z^n$, where here and elsewhere
in this paper $t^\al = t_1^{\al_1} \ldots t_n^{\al_n}$. The Lie
bracket in $G$ is
$$[g_k (\al) , g_r(\be) ] = [g_k, g_r] (\al + \be). $$

\smallskip
{\bf Example 2.}  Let
$R_{n}=\bC[t^{\pm 1}_{1}, \ldots, t^{\pm 1}_{n}]$.
The {\it Witt algebra} (or Cartan type W Lie algebra) is the Lie algebra
$W_{n}=Der(R_{n})=$ Span $\{
t^\a\partial_i|\a  \in \Z^{n}, 1\le i\le n\}$, where
$\partial_i=t_i{\frac{\partial}{\partial t_i}}$.
The bracket in $W_{n}$ is given by
$$[t^\a\partial_i,t^{\be}\partial_{j}]=
\be_i t^{\al + \be} \partial_j - \al_j t^{\al + \be} \partial_i
.$$ Thus $W_n$ is a polynomial Lie algebra.

\smallskip
{\bf Example 3.} The following family of algebras plays an important role in the
representation theory of toroidal Lie algebras. Let $\Omega^1_n$ be the space of
1-forms on a torus: $\Omega^1_n = \sum\limits_{p=1}^n R_n k_p$, where
$k_p = t_p^{-1} dt_p$. We define a 2-parametric family of algebras
$$\V (\mu, \nu) = W_n \oplus \Omega^1_n / d R_n.$$
The distinguished spanning set is
$$\{ d_j (\al), k_j (\al) | \al\in\Z^n, j = 1, \ldots n \},$$
where $d_j (\al) = t^\al \partial_j$, $k_j (\al) = t^\al k_j$. The
Lie bracket in $\V (\mu, \nu)$ is given by
$$[d_i (\al), k_j (\be)] = \be_i k_j(\al+\be) +
\delta_{ij} \sum_{p=1}^n \al_p k_p(\al+\be) ,$$
$$[d_i (\al), d_j (\be)] = \be_i d_j(\ab) - \al_j d_i(\ab) +
(\mu \be_i \al_j + \nu \al_i \be_j) \sum_{p=1}^n \be_p k_p(\al+\be) ,$$
$$ [k_i (\al), k_j (\be)] = 0 .$$
Note that the distinguished spanning set for this polynomial Lie algebra is not a
basis because of the linear dependencies between $k_j (\al)$: \,
$\al_1 k_1(\al) + \ldots + \al_n k_n(\al) = 0$.

\smallskip
Next we will give an example of a Lie algebra which is
exp-polynomial, but not polynomial.
\smallskip

{\bf Example 4.} Consider an associative quantum torus
$$\C_q = \mathop\oplus\limits_{i \in\Z, \al\in\Z^n} \C t_0^{i} t^\al $$
generated by the variables $t_0^\pm, t_1^\pm, \ldots, t_n^\pm$, where
$t_1, \ldots, t_n$ commute and $t_0$ does not commute with $t_1, \ldots, t_n$,
but satisfies the relations:
$ t_j t_0 = q_j t_0 t_j, $
for some $q_1, \ldots, q_n \in \C^*$. The Lie algebra which is obtained from the
associative algebra $C_q$ is an exp-polynomial Lie algebra with the distinguished
set $t_0^i (\al) = t_0^i t^\al$, and the Lie bracket given by
$[t_0^i(\al), t_0^j (\be)] = (q^{j\al} - q^{i\be}) t_0^{i+j}(\ab) .$
Here $q^{j\al} = q_1^{j\al_1} \ldots q_n^{j\al_n}$.

\smallskip
{\bf Definition 1.3.} Let $G=\oplus_{\a\in\Z^n}G_\a$ be an
exp-polynomial Lie algebra. A $G$
module $V=\oplus_{\a\in\Z^n}V_\a$ is called {\it a $\Z^n$-graded
exp-polynomial module} if $V$ has a basis $\{v_j(\a)\}_{j\in
J,\a\in\Z^n}$, and there exists a family of exp-polynomial
functions $h_{k,j}^{s}(\a,\b)$ for $k\in K$, $j,s \in J$ such that
$$g_k(\a)v_j(\b)=\sum_{s\in J} h_{k,j}^{s}(\a,\b) v_s(\a+\b),$$
where $\{g_k(\a)\}_{k\in K}$ is the distinguished spanning set for $G$, and
for each $k,j$ the set $\{s|h_{k,j}^{s}(\a,\b)\ne0\}$ is finite.
The homogeneous components of the $\Z^n$-grading on $V=
\mathop\oplus\limits_{\a\in\Z^n} V_\a$ are given by $V_\a =$ Span
$\{v_j(\a)\}_{j\in J}$.

Note that if $G$ is a  Lie algebra with polynomial multiplication and
 all $h_{k,j}^{s}$ are polynomial functions,
then the module $V$ is actually a polynomial module
defined in [Definition 1.8, BB].

\smallskip
{\bf Example 5.} Let $V_1, \ldots, V_k$ be finite dimensional modules for a
finite dimensional Lie algebra $\gg$. Fix $q_1,\ldots, q_k \in \left( \C^*
\right)^n$. We define on the space $V = R_n \otimes V_1 \otimes \ldots \otimes V_k$
the structure of an exp-polynomial module for the toroidal Lie algebra
$G = R_n \otimes \gg$ by
$$g(\al) v_1 \otimes \ldots \otimes v_k (\be) = \sum_{p = 1}^k q_p^{\al_p}
 v_1 \otimes \ldots (g v_p) \otimes \ldots \otimes v_k (\ab),$$
where $v_1 \otimes \ldots \otimes v_k (\be) = t^\be v_1 \otimes \ldots \otimes v_k$.
It is easy to see that in general $V$ is an exp-polynomial module and not
a polynomial module.

\smallskip
{\bf Example 6.} A tensor module for the Witt algebra $W_n$ is a polynomial module.
Let $V_0$ be a finite-dimensional $gl_n$ module. Then the tensor module
$V = R_n \otimes V_0$
is a module for $W_n$ under the action
$$t^{\al}\partial_j  v(\be) = \be_j v(\ab) + \sum_{p=1}^n \al_p (E_{pj} v) (\ab) ,$$
where $v(\be) = t^\be \otimes v$, and $E_{pj}$ are the elementary matrices in
$gl_n$.
The families $\{ v_i (\be) \}$ with $\{ v_i \}$ being a basis of $V_0$, form
a distinguished basis in $V$, for which the structure constants are polynomials.

\smallskip
We extend  [Definition 1.6, BB] to give the following

{\bf Definition 1.4.}
 Let $G$ be a $\Z^n$-graded exp-polynomial
Lie algebra.  We call this algebra $\Z^n$-{\it extragraded} if
 $G$ has another $\Z$-gradation
 $$G=\mathop\oplus\limits_{i\in\Z }G^{(i)},\eqno(1.2)$$
and the set $K$ is a disjoint union of {\it finite} subsets $K_i$,
$$ K = \mathop\cup\limits_{i\in\Z} K_i,$$
such that
 the elements of the homogeneous spanning set
$\{g_k(\a)\,\,|\,\,k\in K_i,\,\,\a\in\Z^n\}$ are homogeneous
%of degree deg($g_i$)
of degree  $i$ under this new $\Z$-gradation and span $G^{(i)}$.

\smallskip
Many important infinite-dimensional Lie algebras are in fact $\Z^n$-extragraded
exp-polynomial Lie algebras. Here we give some examples.

\smallskip
{\bf Example 8.} We slightly modify Example 1 to get an extragraded algebra. Adding
an extra variable $t_0$, we get an $n+1$ toroidal Lie algebra
$R_{n+1} \otimes \gg$, where $R_{n+1}=\C[t_0^\pm, t_1^\pm, \ldots, t_n^\pm]$.
We consider a $\Z^n$-grading on this algebra by degrees in $t_1, \ldots, t_n$,
and a $\Z$ grading by degree in $t_0$. The distinguished spanning set is
$$t_0^i g_k (\al) = t_0^i t^\al g_k ,$$
with $i \in \Z, \al\in\Z^n$, and $\{ g_k \}$ being a
basis of $\gg$.

\smallskip
The previous example can be generalized in two ways:

{\bf Example 9.} Let $G^{(0)}$ be an exp-polynomial Lie algebra
with $G^{(0)}_\a$ being finite dimensional for all $\a\in\Z^n$. Then
$ G = \C[t_0, t_0^{-1}] \otimes G^{(0)}$
is an extragraded exp-polynomial Lie algebra.

\smallskip
{\bf Example 10.} Let $\gg$ be a $\Z$-graded Lie algebra with finite-dimensional
homogeneous components. Then
$ G = \gg \otimes R_n $
is an extragraded exp-polynomial Lie algebra.

\smallskip
{\bf Remark.}  In the definition of an exp-polynomial Lie algebra we can
relax the requirement that each family
$\{ g_k(\al) \}$ be defined for all $\al\in\Z^n$. We may instead require that
$\{ g_k(\al) \}$ is defined for $\al$ in some sublattice $L_k \subset \Z^n$.
Of course in this case we should have a restriction so
that the expression  in the right hand side of (1.1) is well defined.

\smallskip
{\bf Example 11.} With this relaxed definition we may consider the following as an
extragraded exp-polynomial Lie algebra:
$$ G = \gg \otimes R_{n+1} \oplus Der \, \C[t_0, t_0^{-1}] $$
for a finite dimensional Lie algebra $\gg$ with the distinguished spanning set
$\{ t_0^i g_k (\al),
\break
t_0^i \partial_0  | i \in\Z, \al\in\Z^n, g_k\in B \}$.
Here the
sublattice that corresponds to the elements $t_0^i \partial_0$ is $L = (0)$.

\smallskip
{\bf Example 12.}
Introducing an extra variable $t_0$ we can construct an extragraded version of the
Lie algebra from Example 3:
$$\V(\mu,\nu) = Der R_{n+1} \oplus \Omega^1_{n+1} / dR_{n+1} .$$
The $\Z$-grading is by degree in $t_0$, and the distinguished
spanning families are $t_0^i d_j (\al) = t_0^i t^\al \partial_j$
and $t_0^i k_j (\al) = t_0^i t^\al k_j$.

\smallskip
{\bf Example 13.} The exp-polynomial Lie algebra from Example 4 is actually
extragraded. The $\Z$-grading on it is by degree in $t_0$.

\smallskip
{\bf Example 14.}  Let
$R_{n+1}=\bC[t^{\pm 1}_{0},t^{\pm 1}_{1}, \ldots, t^{\pm 1}_{n}]$,
$W_{n+1}=Der(R_{n+1})=$ \break Span
$\{t_0^i t^\a\partial_k, t_0^i t^\a\partial_0|i\in\Z, \a \in \Z^{n}, 1\le k\le n\}$.
% where
% $t^\a=t_0^{\a_0}t_1^{\a_1}...t_n^{\a_n}$,  $\partial_i=t_i{\frac{\partial}
% {\partial t_i}}$.
The Lie bracket in $W_{n+1}$ is given by
$$[t_0^i t^\a\partial_k, t_0^j t^{\b}\partial_{s}] =
t_0^{i+j} t^{\a+\b} (\b_{k}\partial_{s}-\a_{s}\partial_{k}), \, \, 1\le k,s \le n,$$
$$[t_0^i t^\a\partial_0, t_0^j t^{\b}\partial_{s}] =
t_0^{i+j} t^{\a+\b} (j \partial_{s}-\a_{s}\partial_{0}), \, \, 1\le s \le n,$$
$$[t_0^i t^\a\partial_0, t_0^j t^{\b}\partial_{0}] =
(j-i) t_0^{i+j} t^{\a+\b} \partial_{0}.$$
It can be easily seen that $W_{n+1}$ is a $\Z^n$-extragraded  Lie algebra.

\smallskip
Cartan type S Lie algebras $S_{n+1}$ are also $\Z^n$-extragraded
polynomial Lie algebras. However we do not know whether Cartan
type H or K Lie algebras are $\Z^n$-extragraded polynomial Lie
algebras.

%[Kaiming, could you please expand on the above remark? Other types?]

\smallskip
{\bf Example 15.}  The {\it Virasoro-like algebra} $L$ over $\bC$
is the Lie algebra with a $\bC$-basis $\bigl\{L_{x}\bigm|x\in \Z^2
\bigr\}$ and subject to the following commutator relations
$$[L_x, L_y]=\det
\left(\matrix y\crcr x\endmatrix\right)
L_{x+y},\qquad\forall\,\,x,y\in
\Z^2,
$$
where $x=(x^{(1)},x^{(2)}), y=(y^{(1)},y^{(2)})$, $\left(\matrix
y\crcr x\endmatrix\right) = \left(\matrix y^{(1)}&y^{(2)}\crcr
x^{(1)}&x^{(2)}\endmatrix\right)$.
It can be easily seen that $L$ is a $\Z^1$-extragraded  Lie algebra.

\smallskip
From now on we assume that $G$ is a $\Z^n$-extragraded  Lie algebra with gradations (1.1) and (1.2), i.e.,
$G=\oplus_{i\in\Z,\,\,\,\a\in\Z^{n}}G^{(i)}_\a$ is a $\Z^{n+1}$-graded Lie algebra
which has a homogeneous spanning set
$\{g^{(i)}_k(\a)\,\,\,|\,\,\,{k\in K_i,\,\,(i,\a)\in\Z^{n+1}}\}$ with
 $g^{(i)}_k(\a)\in G^{(i)}_\a$,
and there exists a family of exp-polynomial functions $\{f_{k,m,i,j}^s(\a,\b)\}$ in
the $2n$ variables $\a_p,\b_p$ where
$k \in K_i, m\in K_j, s\in K_{i+j}$ and where for each $k,m,i,j$ the set $\{s|f_{k,m,i,j}^s(\a,\b)\ne0\}$ is
finite, such that
$$[g^{(i)}_k(\a),g^{(j)}_m(\b)]=\sum_{s\in
K_{i+j}}f_{k,m,i,j}^s(\a,\b)g^{(i+j)}_s(\a+\b),\; \f \a,\b\in\Z^n.\eqno(1.3)$$

Let $G^+=\mathop\oplus\limits_{i\ge1}G^{(i)}$,
$G^-=\mathop\oplus\limits_{i\le-1}G^{(i)}$.
Then we have the decomposition
$$G=G^-\oplus G^{(0)} \oplus G^+.\eqno(1.4)$$
\smallskip
Note that $G^{(0)}$ is a $\Z^{n}$-graded exp-polynomial Lie algebra.

\smallskip
Following the construction in [BB], we now introduce our $\Z^{n+1}$-graded module over
$\Z^{n+1}$-graded Lie algebra $G$.

\smallskip
Assume $V=\oplus_{\a\in\Z^n}V_\a$ is a  $\Z^n$-graded $G^{(0)}$ module with
exp-polynomial action as defined  in Definition 1.3.
We can define the action of $G^+$ on $V$ by $G^+V=0$ and then consider the induced module
$$\tilde M(V)={\roman{Ind}}^G_{G^{(0)}+G^+}V\simeq U(G^-)\otimes_\bC V. \eqno{(1.5)}$$
It is clear that $\tilde M(V)$ is a $\Z^{n+1}$-graded module over $G$ and
$$\tilde M(V)=\mathop\oplus\limits_{i\le0, \a\in\Z^{n}}\tilde
M(V)^{(i)}_{\a},\eqno(1.6)$$
where $\tilde M(V)^{(i)}_{\a}$ is naturally defined, for example,
$\tilde M(V)^{(0)}_{\a}=V_{\a}$.
In general, the homogeneous components $\tilde M(V)^{(i)}_{\a}$ with $i < 0$ are
infinite-dimensional.

\smallskip
It is easy to see  that $\tilde M(V)$ has a unique maximal proper $\Z^{n+1}$-graded
submodule $\tilde M^{rad}$ which intersects trivially with $V$.
Let
$$M(V)=\tilde M(V)/\tilde M^{rad}.\eqno{(1.7)}$$
Then we have the induced $\Z^{n+1}$-gradation
$$M(V)=\mathop\oplus\limits_{i\le0, \a\in\Z^{n}}M(V)^{(i)}_{\a}.\eqno(1.8)$$

\smallskip
The main result of this paper is the following theorem.

 \medskip
{\bf Theorem 1.5.}   {\it Assume that $G$ is a $\Z^n$-extragraded
Lie algebra with grading (1.3), $V=\oplus_{\a\in\Z^n}V_\a$ is a
$\Z^n$-graded exp-polynomial $G^{(0)}$ module as defined in
Definition 1.3 with $J$ being finite. Then the $\Z^{n+1}$-graded
$G$ module $M(V)$ defined in (1.7) has finite-dimensional
homogeneous spaces, i.e., $\dim M(V)^{(i)}_{\a} <\infty,$ for all
$ i\in \Z,\, \break \a\in \Z^n.$}
\smallskip
If $G^{(0)}$ is a polynomial Lie
algebra and $V$ is a polynomial $\Z^n$-graded $G^{(0)}$ module,
the claim of Theorem 1.5 was proved in [Theorem 1.12, BB].

\smallskip
{\bf Example 16.}  Let us consider an $n+1$-toroidal Lie algebra
$G=R_{n+1}\otimes \gg$ as defined in Example 8. Its component
$G^{(0)}$ is an $n$-toroidal Lie algebra. We consider a module $V =
R_n \otimes V_1 \otimes \ldots \otimes V_k$ for $G^{(0)}$ which was
described in Example 5. By the above theorem, the $G$ module
$M(V)$ has finite-dimensional homogeneous components. Such modules
were studied by Chari [C] and Rao [E].

\smallskip
{\bf Example 17.}  Let us consider the Virasoro-like algebra $L$
as defined in Example 15.  Let
$L^{(i)}=\oplus_{k\in\Z}\bC L_{i,k}$.
Fix an exp-polynomial function $f(k)$.
Then $V=\oplus_{j\in\Z}\bC v_j$ becomes
an exp-polynomial $L^{(0)}$ module via
$$L_{0,k}v_j= f(k) v_{j+k}.$$
Then $M(V)$ has  finite dimensional homogeneous spaces.

\smallskip
{\bf Example 18.} Consider a Witt algebra $W_{n+1} = Der \, \C
[t_0^\pm, t_1^\pm, \ldots, t_n^\pm]$, which we view as an
extragraded polynomial Lie algebra with $\Z$ grading given by
degree in $t_0$. The zero component with respect to this grading
is $ W_n \oplus R_n \partial_0 .$ Consider a tensor module $V =
R_n \otimes V_0$ for $W_n$ as discussed in Example 6. We let $R_n
\partial_0$ act upon it by shifts $ t^\al \partial_0 \cdot t^\be
\otimes v = d \, t^\ab \otimes v $ for some fixed constant $d \in
\C$. By Theorem 1.5 the module $M(V)$ is a weight module with
finite dimensional weight spaces.

\smallskip
{\bf Example 19.} Let $\V(\mu,\nu)$ be the extragraded Lie algebra from
Example 12. Its zero component with respect to $\Z$ grading is
$$ W_n \oplus R_n \partial_0 \oplus \left( \sum_{p=0}^n R_n k_p / d R_n \right) .$$
Consider a tensor module $V$ from the previous example on which we define the action
of 1-forms of degree zero (in $t_0$) as follows:
$$ t^\al k_p \cdot t^\be \otimes v = 0, \; p = 1, \ldots, n,$$
$$ t^\al k_0 \cdot t^\be \otimes v = c t^\ab \otimes v,
\hbox{\rm \ for some \ } c \in\C .$$ Again by Theorem 1.5 the
module $M(V)$ is a weight module with finite dimensional weight
spaces. These modules and their irreducible quotients were studied
in [EM], [L], [BB], [B1], [B2].

\smallskip
{\bf Example 20.}
Let $G$ be the quantum torus Lie algebra from Example 4. We consider a module
$V = \C[t_1^\pm, \ldots t_n^\pm]$
for $G^{(0)} = \C[t_1^\pm, \ldots t_n^\pm]$ with the action defined as follows:
$t^\al \cdot t^\be = f(\al,\be) t^\ab ,$
where $f(\al,\be)$ is some fixed exp-polynomial function. By Theorem 1.5
the module $M(V)$ has finite-dimensional homogeneous components.

\smallskip
{\bf Example 21.}
Let $\gg$ be a simple finite-dimensional Lie algebra with the triangular decomposition
$\gg = \gg_- \oplus \hh \oplus \gg_+$. We consider a finite $\Z$-grading on $\gg$
compatible with this decomposition. As explained in Example 10 above, the algebra
$G = R_n \otimes \gg$ is an extragraded polynomial Lie algebra. Its zero component
is the abelian algebra $G^{(0)} = R_n \otimes \hh$. Consider a $G^{(0)}$ module
$V = \C [t_1^\pm, \ldots t_n^\pm]$ with the action
$ t^\al h_k \cdot t^\be = f_k(\al,\be) t^\ab ,$
where $\{ h_k \}$ forms a basis of $\hh$ and $\{ f_k (\al,\be) \}$ are some fixed
exp-polynomial functions. Applying again Theorem 1.5 we conclude that
$M(V)$ is a weight module with finite dimensional weight spaces.

\medskip
We would like to discuss now a finite-dimensional version of the exp-polynomial
modules.

\

{\bf Definition 1.6.} Let $G=\oplus_{\a\in\Z^n}G_\a$ be an
exp-polynomial Lie algebra as defined in Definition 1.2. A
$G$-module $V$ is called {\it a finite-dimensional exp-polynomial
module} if $V$ has a finite basis $\{v_j\}_{j\in J}$, and there
exists a family of exp-polynomial functions $h_{k,j}^{s}(\a)$ for
$k\in K$, $j,s \in J$ such that
$$g_k(\a)v_j = \sum_{s\in J} h_{k,j}^{s}(\a) v_s,$$
where $\{g_k(\a)\}_{k\in K, \a\in\Z^n}$ is the distinguished spanning set for $G$.

\

Let $G$ be an extragraded exp-polynomial Lie algebra and let $V$ be a
finite-dimensional exp-polynomial module for $G^{(0)}$. Just as in (1.5) and (1.7) we
define a $\Z$-graded $G$ module $\tilde M(V)$ and its $\Z$-graded factor-module
$M(V) = \tilde M (V)/ \tilde M^{rad}$.

\medskip
{\bf Theorem 1.7.}   {\it Let $G$ be a $\Z^n$-extragraded
Lie algebra with grading (1.3), $V$ be a
finite-dimensional exp-polynomial $G^{(0)}$ module.
Then the
$\Z$-graded $G$ module $M(V)$ has
finite dimensional homogeneous spaces, i.e., $\dim M(V)^{(i)} <\infty,$ for all
$ i\in \Z.$}

\medskip
Let us now give some examples of the applications of this theorem.

{\bf Example 22.}
 Let $G$ be an $n+1$-toroidal Lie algebra
$G=R_{n+1}\otimes \gg$ as defined in Example 8. Its component
$G^{(0)}$ is an $n$-toroidal Lie algebra.
Let $V_1,\ldots, V_k$ be finite-dimensional modules for $\gg$,
and let $q_1,\ldots, q_k \in \left( \C^* \right)^n$.
We define a structure of a $G^{(0)}$ module on the finite-dimensional
space $V = V_1 \otimes \ldots \otimes V_k$ in the following way:
$$g(\al) v_1 \otimes \ldots \otimes v_k = \sum_{p = 1}^k q_p^{\al_p}
 v_1 \otimes \ldots (g v_p) \otimes \ldots \otimes v_k .$$
It is easy to see that $V$ is a finite-dimensional exp-polynomial module
for $G^{(0)}$. The induced module $\tilde M(V)$ will have infinite-dimensional
homogeneous components, however, by Theorem 1.7, the homogeneous components
of its factor $M(V)$ are finite-dimensional.
% See [C] for more details on these modules.

The next two examples are modifications of Examples 20 and 21.

\smallskip
{\bf Example 23.}
Let $G$ be the quantum torus Lie algebra from Example 4. We consider a one-dimensional
module
$\C$ for $G^{(0)} = \C[t_1^\pm, \ldots t_n^\pm]$ with the action defined as follows:
$t^\al \cdot 1 = f(\al) 1 ,$
where $f(\al)$ is some fixed exp-polynomial function (highest weight). By Theorem 1.7
the module $M(\C)$ has finite-dimensional homogeneous components.

\smallskip
{\bf Example 24.}
Let $\gg$ be a simple finite-dimensional Lie algebra with the triangular decomposition
$\gg = \gg_- \oplus \hh \oplus \gg_+$. As explained in Example 21 above, the algebra
$G = R_n \otimes \gg$ is an extragraded polynomial Lie algebra. Its zero component
is the abelian algebra $G^{(0)} = R_n \otimes \hh$.
Let $\{ h_1, \ldots, h_\ell \}$ be a basis of $\hh$.
Fix exp-polynomial functions $f_1(\al), \ldots, f_\ell(\al)$ (highest weight) and
consider a one-dimensional  $G^{(0)}$ module $\C$ with the action
$ t^\al h_k \cdot 1 = f_k(\al) 1$.
Applying again Theorem 1.7 we conclude that
$M(\C)$ is a $\Z$-graded module with finite dimensional homogeneous components.

\smallskip\bigskip
\subhead 2. \ \ Proof of the main theorems 
\endsubhead
\medskip

In this section we shall prove Theorems 1.5 and 1.7.
%Our proof will be quite
%parallel to the proof of Theorem 1.12 from [BB], generalizing the
%polynomial case to the exp-polynomial case.
The key step in the
proof of Theorem 1.12 in [BB] was the Vandermonde determinant
argument. Here we would need a generalization of the Vandermonde
determinant formula for the exp-polynomial functions.

\medskip

{\bf Lemma 2.1.} {\it Let $a_1,\ldots,a_m$ be elements of a field,
$s_1,s_2,\ldots, s_m \in \N$ with  $s_1 + \ldots + s_m = s$. Consider the
following sequence of $s$ exp-polynomial functions in one integer variable:
$f_1(n) = a_1^n, f_2(n) = n a_1^n, \ldots, f_{s_1}(n) = n^{s_1-1} a_1^n$,
$f_{s_1+1} (n) = a_2^n, \ldots, f_{s_1+s_2}(n) = n^{s_2-1} a_2^n$,
$\ldots, f_s (n) = n^{s_m -1} a_m^n$.
Let $V=(v_{pk})$ be the  square $s\times s$ matrix  where $v_{pk} = f_k (p-1)$,
$p,k = 1, \ldots s$. Then
$$ \det(V) = \mathop\Prod\limits_{j=1}^m (s_j - 1)!! \; a_j^{\frac{s_j (s_j-1)}{2}} \;
\mathop\Prod\limits_{1 \leq i < j \leq m} (a_j - a_i)^{s_i s_j} .\eqno(2.1)$$}

Here we use the notation $m!! = m! \times (m-1)! \times
\ldots \times 2! \times 1!$ with the convention $0!!=1$.

 \vs
{\it Proof.} This elementary lemma may be proved by
induction on $s$ using elementary row and column transformations
of the matrix. We will give here just an outline of the proof.
The basis of induction case $s=1$ is trivial.
To establish the inductive step, one begins by applying
elementary row operations to the matrix, the same as in the proof of the ordinary
Vandermonde determinant formula: subtract from the last row the preceding row number
$s-1$ multiplied by $a_1$, then subtract from row $s-1$ the preceding row $s-2$
multiplied by $a_1$, and so on. This will produce a matrix that has $1$ as the top
entry of the first column, and the rest of the entries in the first column being
zeros.

Next we expand this determinant along the first column, which will yield an
\break
$(s-1) \times (s-1)$ matrix with the same determinant. Finally, applying elementary column
operations, it is possible to bring this $(s-1)\times (s-1)$ matrix to the form
corresponding to the sequence of functions
$a_1 a_1^n, \, 2 a_1 n a_1^n, \, 3 a_1 n^2 a_1^n, \ldots,
\break
(s_1-1) a_1 n^{s_1 -2}
a_1^n, \,
(a_2-a_1) a_2^n, \, (a_2-a_1) n a_2^n, \ldots, (a_2-a_1) n^{s_2-1} a_2^n, \,
(a_3-a_1) a_3^n, \ldots,
\break
(a_m-a_1) n^{(s_m-1)} a_m^n$.

Pulling out the factors $a_1, 2 a_1, 3 a_1, \ldots, (s_1-1) a_1$ from the first $s_1-1$
columns, and $(a_j-a_1)$ from the remaining columns, we bring the matrix to the extended
Vandermonde form of rank $s-1$. To establish the inductive step, we only need to multiply
the expression for the extended Vandermonde determinant of rank $s-1$
given by the induction assumption by the factors we pulled out of the columns.
\qed

\vs

{\bf Corollary 2.2.} {\it In the notations of the previous lemma, let
$a_1,\ldots,a_m$ be
distinct non-zero elements
of a field of characteristic greater or equal to the maximum of $s_1,
\ldots s_m$, or of characteristic $0$. Then the set of exp-polynomial
functions $f_1(n) = a_1^n, \ldots, f_s(n) = n^{s_m-1} a_m^n$ is linearly
independent.}
\vs
{\it Proof}. From the determinant formula (2.1) we see that the vectors of values
of the functions $\left(f_j(0), f_j(1), \ldots, f_j(s-1)\right)$, $j=1, \ldots, s$, are
linearly independent. \qed

\vs

{\bf Corollary 2.3.} {\it Let the exp-polynomial functions $f_1(n),\ldots f_s(n)$
satisfy the conditions of Lemma 2.1 and Corollary 2.2. Let $\{ c_k \}$ be
a sequence with only finitely many non-zero terms. The sequence $\{ c_k \}$
satisfies an infinite system of linear equations
$$\sum_k \left( \sum_{j=1}^s d_{kj} f_j(n) \right) c_k = 0
\hbox{\rm \ \ for all \ } n \in\Z  \eqno(2.2) $$
if and only if it satisfies a finite system
$$ \sum_k d_{kj} c_k = 0 \hbox{\rm \ \ for all \ } j = 1, \ldots, s. \eqno(2.3)$$}
\vs
{\it Proof}. Since
$$\sum_k \left( \sum_{j=1}^s d_{kj} f_j(n) \right) c_k =
\sum_{j=1}^s f_j(n) \sum_k d_{kj} c_k,$$
we see that (2.3) implies (2.2).
\vs
Suppose now that (2.2) holds. Evaluating (2.2) at $n=0,\ldots,s-1$,
we get that
$$\sum_{j=1}^s f_j(n) \sum_k d_{kj} c_k = 0 \hbox{\rm \ for \ }
n=0,\ldots,s-1. \eqno(2.4)$$
Now $s$ equations of (2.4) are linear combinations of $s$ equations of
(2.3). The change of basis matrix from (2.3) to (2.4) is an extended
Vandermonde matrix and is invertible,
since by Lemma 2.1 it has a non-zero determinant. Hence, the equations
of (2.3) are linear combinations of equations of (2.4) and thus (2.3)
holds.\qed
\vs
The last Corollary also admits a straightforward multi-variable
generalization.
\vs
{\bf Corollary 2.4.} {\it Let $f_1(n_1,\ldots,n_r), \ldots f_s(n_1,\ldots, n_r)$
be a set of distinct exp-polynomial functions of the form
$f_j (n_1,\ldots,n_r) = n_1^{p_1} \ldots n_r^{p_r} b_1^{n_1} \ldots
b_r^{n_r}$, such that exponents $p_i$'s are less than the characteristic
of the field, in case if the field is of finite characteristic.
 Let $\{ c_\b \}_{\b\in\Z^m}$ be
a set with only finitely many non-zero terms. The set $\{ c_\b \}$
satisfies an infinite system of linear equations
$$\sum_{\b\in\Z^m} \left( \sum_{j=1}^s d_{\b j} f_j(n_1, \ldots, n_r) \right) c_\b = 0
\hbox{\rm \ \ for all \ } (n_1, \ldots, n_r) \in\Z^r  \eqno(2.5) $$
if and only if it satisfies a finite system
$$ \sum_{\b\in\Z^m} d_{\b j} c_\b = 0 \hbox{\rm \ \ for all \ } j = 1, \ldots, s.
\eqno(2.6)$$}

\vs
Now we are ready to give a proof of Theorem 1.5.

\medskip
{\it Proof of Theorem 1.5.}  For the $\Z^n$-graded $G^{(0)}$ module
$V=\oplus_{\a\in\Z^n}V_\a$ as defined in Definition 1.3, we have
used $g_i(\a)$ to denote $g^{(0)}_i(\a)$. We stress here that both
$\{k\,\,|\,\,g^{(i)}_k(\a)\ne0\}$ for any fixed
$(i,\a)\in\Z^{n+1}$ and $J$ in Definition 1.3  are finite sets.
% We shall simply denote $M(V)$ by $M$.
The proof of the theorem will amount to proving two claims.
\medskip
{\bf Claim 1}.
{\it Let us fix $i_1,i_2,\ldots,i_{s}\in\N$,
$k_1,k_2,\ldots,k_{s}$, with $k_p \in K_{-i_p}$, $k\in J$ and
$\a\in\Z^n$. There exist exp-polynomial functions
$f_1(\b),\ldots,f_d(\b)$ in $ns$ variables $\beta$ such that a
linear combination
$$\sum_{\b=(\b_1,\ldots,\b_s)\in\Z^{ns}}b_\b \, g_{k_1}^{(-i_1)}(\b_1)
g_{k_2}^{(-i_2)}(\b_2)\ldots g_{k_s}^{(-i_s)}(\b_s)
v_{k }(\a-\b_1-\ldots-\b_s) \eqno(2.7)$$
belongs to $\tilde M^{rad}$ if and only if the set $\{ b_\b \}_{\b\in\Z^{ns}}$
(with finitely many non-zero elements)
satisfies a finite system of linear equations:
$$\sum_{\b\in\Z^{ns}}b_\b f_p (\b) = 0, \hbox{\rm \ for \ }
p = 1, \ldots, d.$$
}

{\it Proof of Claim 1.} Denote the sum in (2.7) by $x$. We have that
$x\in\tilde M^{rad}$ if and only if
$$ g_{m_1}^{(j_1)}(\g_1) \ldots g_{m_r}^{(j_r)}(\g_r) x = 0$$
for all $j_1,\ldots,j_r\in\N$ with $j_1+\ldots+j_r = i_1 + \ldots + i_s$,
$m_1\in K_{j_1},\ldots,m_r\in K_{j_r}$, $\g_1, \ldots, \g_r \in \Z^n$.

From the Poincar\'e-Birkhoff-Witt argument and the fact that the Lie algebra $G$
and the module $V$ are exp-polynomial, it follows that
$$ g_{m_1}^{(j_1)}(\g_1) \ldots g_{m_r}^{(j_r)}(\g_r) x =
\sum_{\b\in\Z^{ns}} b_\b \sum_{\ell\in J} h_\ell (\b,\g)
v_\ell (\alpha + \g_1 + \ldots \g_r) .$$
The functions $h_\ell$ are exp-polynomial in $\b = (\b_1,\ldots\b_s)$
and $\g = (\g_1,\ldots,\g_r)$, and depend on $j_1,\ldots,j_r; m_1, \ldots, m_r$.
We get that $x\in\tilde M^{rad}$ if and only if
$\{ b_\b \}$ satisfies the system of equations
$$\sum_{\b\in\Z^{ns}} b_\b \sum_{\ell\in J} h_\ell (\b,\g) = 0 \eqno{(2.8)}$$
for all $\g\in \Z^{nr}, \ell\in J, j_1,\ldots,j_r\in\N$
with $j_1+\ldots+j_r = i_1 + \ldots + i_s$,
$m_1\in K_{j_1},\ldots,m_r\in K_{j_r}$.
Since $j_1,\ldots,j_r\in\N$ are bounded by the condition
$j_1+\ldots+j_r = i_1 + \ldots + i_s$
and $m_1, \ldots, m_r$ belong to finite sets, we conclude that only finitely many
functions $h_\ell (\b,\g)$ appear in the system (2.8). Nonetheless the system
(2.8) has infinitely many equations because $\g$ has an infinite range.
Our goal is to reduce (2.8) to a system with finitely many equations.

A exp-polynomial function $h_\ell (\b,\g)$ could be expanded in $\g$:
$$h_\ell (\b,\g) = \sum_p f_p(\b) a_p^\g \g^{\d_p}, $$
where $a_p \in \C^{nr}, \d_p \in \Z^{nr}$ and the summation in $p$ is finite.
By Corollary 2.4, the system (2.8) is equivalent to a finite system of linear equations:
$$\sum_{\b\in\Z^{ns}} b_\b f_p(\b) = 0,$$
with a finite number (say, $d$) of exp-polynomial functions $f_p(\b)$. This
establishes Claim 1.

\

{\bf Claim 2.}
{\it In the notations of Claim 1, the set
$$ \left\{ g_{k_1}^{(-i_1)}(\b_1) \ldots g_{k_s}^{(-i_s)}(\b_s)
v_{k }(\a-\b_1-\ldots-\b_s) \right\}_{\b\in\Z^{ns}} $$
spans a subspace of dimension less or equal to $d$ in $M_\a^{(-i_1-\ldots-i_s)}(V)$.
}

{\it Proof of Claim 2.} To prove this claim we need to show that any $d+1$ vectors
in this set are linearly dependent. Let $B$ be a subset in $\Z^{ns}$ of size $d+1$.
The homogeneous system of $d$ linear equations
$$\sum_{\b\in B} b_\b f_p (\b) = 0$$
in $d+1$ variables $\{ b_\b \}_{\b\in B}$ has a non-trivial solution. By Claim 1,
the set
$$\left\{ g_{k_1}^{(-i_1)}(\b_1) \ldots g_{k_s}^{(-i_s)}(\b_s)
v_{k }(\a-\b_1-\ldots-\b_s) \right\}_{\b\in B} $$
is linearly dependent modulo $\tilde M^{rad}$. Claim 2 is now proved.

\

Theorem 1.5 now follows immediately from Claim 2 since the space
$M_\a^{(-i)}$ is spanned by the elements
$$\left\{ g_{k_1}^{(-i_1)}(\b_1) \ldots g_{k_s}^{(-i_s)}(\b_s)
v_{k }(\a-\b_1-\ldots-\b_s) \right\}$$
with $i_1+\ldots+i_s = i; i_1,\ldots,i_s \in\N$. Thus there are finitely many
possibilities for $(i_1,\ldots,i_s)$, and the indices $k_1,\ldots,k_s$ and $k$
run over finite sets $K_{-i_1},\ldots,K_{-i_s}$ and $J$.
\qed

{\it Proof of Theorem 1.7.} The proof of Theorem 1.7 may be derived from the proof of
Theorem 1.5 by forgetting $\Z^n$-grading on the module $V$ and replacing
expressions $v_{k }(\a-\b_1-\ldots-\b_s)$ etc., simply with $v_k$.
\qed

\smallskip\bigskip
\subhead 3. \ \ Weight modules over generalized Virasoro algebras
\endsubhead
\medskip

The {\it Virasoro algebra} $\Vir:=\Vir[\Z]$
over $\bC$ is a Lie algebra with a $\bC$-basis
$\bigl\{c,d_{i}\bigm|i\in\Z\bigr\}$ and subject to the following
commutator relations $$\eqalign{ &[c, d_i]=0,\cr &[d_i,
d_j]=(j-i)d_{i+j}+\delta_{i,-j}\frac{i^{3}-i}{12}c,
\qquad\f \,\,i,j\in\Z.\cr}
$$
The structure and representation theory for the Virasoro algebra
has been well developed. For details, we refer the readers to [M]
and the book [KR] and the references therein.
%The centerless Virasoro
%algebra is essentially a Witt algebra, and Witt algebras in
%positive characteristic were studied by many authors, for
%instance, Zassenhaus [Z], Strade [St], Ree [R] and Wilson [W].

\smallskip
{\it Generalized Virasoro algebras} $\Vir[M]$ in characteristic $0$ were
introduced by Patera and Zassenhaus in [PZ], which are Lie
algebras obtained from the Virasoro algebra
by replacing the index group $\Z$ with an
arbitrary nonzero subgroup $M$ of the base field $\bC$, which is also the
1-dimensional universal central extension of some generalized Witt
algebras [DZ]. Representations for Generalized Virasoro algebras
$\Vir[M]$ have been studied in several references. In [SZ], all
weight modules with weight multiplicity $1$ were determined, which
are some so-called intermediate series modules. In [Ma], it was
proved that all irreducible Harish-Chandra modules  over the
generalized Virasoro algebra $\Vir[\bQ]$ over $\bC$ are
intermediate series modules (where $\bQ$ is the rational number
field).
In [HWZ], the
irreducibility of Verma modules over the
generalized Virasoro algebra $\Vir[M]$ over $\bC$ was
completely determined.

%\smallskip
%In [S], Su claimed that, if $M\simeq \Z^n$ with $n>1$, then any
%irreducible Harish-Chandra module   (i.e., irreducible weight
%module with finite weight multiplicity) over  Vir$[M]$ is a module
%of the intermediate series (i.e., all its weight multiplicities
%are $1$), thus gave the complete classification for Harish-Chandra
%modules over the higher rank Virasoro  algebras. Unfortunately,
%this is not correct.

\smallskip
In this section we assume that
$M= \Z\oplus M_0\subset\bC$ where $M_0$ is a nonzero subgroup of $\bC$.
We simply denote $L=\Vir[M]/\bC c$. For any $i\in\Z$, we denote
$$L_i=\oplus_{a\in M_0}\bC d_{i+a},$$
$$L_+=\oplus_{i\in\Z_+} L_i,\,\, L_-=\oplus_{i<0} L_i.$$
Then $L_0\simeq \Vir[M_0]/\bC c$.

For any $\a,\b\in\bC$, we have the $L_0$ module
$V(\a,\b)=\oplus_{a\in M_0}\bC v_a$ subject to the actions
$$d_av_b=(b+\a+a\b)v_{a+b},\f a,b\in M_0.\eqno(3.1)$$
We extend  $L_0$ module structure on
$V(\a,\b)$ to  $L_++L_0$ module structure by defining
$$L_+V(\a,\b)=0.\eqno(3.2)$$
Then we obtain the induced $L$ module
$$\tilde V=\tilde V(\a,\b)
=\Ind_{L_++L_0}^{L}V(\a,\b)=U(L)\otimes_{U(L_++L_0)}V(\a,\b),$$
 where $U(L)$ is the universal enveloping algebra of the Lie algebra $L$.
 It is clear that, as vector spaces,
 $ \tilde V\simeq U(L_-)\otimes_{\bC}V(\a,\b)$. The   module $\tilde V$
 has a unique maximal proper submodule $J$ trivially
 intersecting with $V(\a,\b)$. Then we obtain the quotient
 module
 $$\bar V=\bar V(\a,\b)=\tilde V / J.\eqno(3.3)$$
 It is clear that $\bar V$ is uniquely determined by the constants $\a,\b$,
  and
  $\bar V=\break\oplus_{i\in\Z_+}\bar V_{-i+\a+M_0}$ where
 $$\bar  V_{-i+\a+M_0}=\oplus_{a\in M_0}\bar V_{-i+\a+a},\,\,\,
 \bar V_{-i+\a+a}=\{v\in \bar  V\,\,|\,\,d_0v=(-i+\a+a)v\}.\eqno(3.4)$$
We can similarly define $\tilde  V_{i+\a+M_0}$ and $\tilde V_{-i+\a+a}$.
It is easy to see that
$$\dim \tilde V_{-i+\a+a}=\infty,\f i\in \N,\, a\in M_0.$$
It will be different for $\bar V_{-i+\a+a}$.

Note that if $M_0\simeq \Z^n$, then $L$ is a $\Z^n$-extragraded  Lie algebra,
$V(\a,\b)$ is a $L$ module with
polynomial action, and $\bar V(\a,\b)=M(V(\a,\b))$ as defined in Sect.2.
In this case, from Theorem 1.5, we know that
$\bar V(\a,\b)$ has finite dimensional weight spaces , i.e.,
$\dim \bar V_{-i+\a+a} <\infty$ for all $ i\in \N,\, a\in M_0.$
More generally, we have

 \medskip
{\bf Theorem 3.1}   {\it The module $\bar V(\a,\b)$ defined in
(3.3) over $\Vir[M]$ has finite dimensional weight spaces, more
precisely, $\dim \bar V_{-i+\a+a}\le (2i+1)!!$ for all $ i\in
\N,\, a\in M_0.$}

\medskip
{\it Proof.}  Since $L_-$ is generated by $L_{-1}$ ,
and $L_+$ is generated by $L_{1}$,  we deduce  that
$$L_{-1}\bar V_{-i+\a+M_0}=\bar V_{-i-1+\a+M_0},\f i\in\Z_+,$$
and, if $v\in \bar V_{-i+\a+M_0}$ where $i\in\N$, satisfies $L_1v=0$ then $v=0$.
We also know that, for any $ n\in \N,\, a\in M_0$,
$$\bar V_{-n+\a+a}=span\{d_{-1+a_n}d_{-1+a_{n-1}}...d_{-1+a_1}v_{a_0}
\,\,|\,\,a_i\in M_0,$$
$$\hskip 6cm \,\,\,with\,\,\,a_0+a_1+...+a_n=a\}.$$

\medskip
{\bf Claim 1}. {\it For any $n\in\Z_+$, and $b_i\in\bC$, $\a_i\in M_0$,
(only finitely many $b_i$ are not zero),

(a) if
$$\sum_{i\in\Z}\a_i^kb_i=0,\f\,\,0\le k\le 2n+2,\eqno(3.5)$$
then $$\sum_{i\in\Z}b_id_{-1+\a_i}d_{-1+a_{n}}...d_{-1+a_1}v_{a_0-\a_i}=0;
\eqno(3.6)$$

(b)  if
$$\sum_{i\in\Z}\a_i^kb_i=0,\f\,\,0\le k\le 2n+1,\eqno(3.5')$$ then
$$\sum_{i\in\Z}b_id_{\a_i}d_{-1+a_{n}}...d_{-1+a_1}v_{a_0-\a_i}=0,\eqno(3.7)$$
for all $a_0,a_1,...,a_n\in M_0$, where all the sums are finite.}
\medskip

Now we first consider $n=0$. Suppose $(3.5')$ holds for $n=0$. We
have
$$\sum_{i\in\Z}b_id_{\a_i}v_{a_0-\a_i}=\sum_{i\in\Z}b_i(a_0-\a_i+\a+\b\a_i )
v_{a_0}=0,$$
which is (3.7) for $n=0$. Now suppose (3.5) holds for $n=0$.
For any $\g\in M_0$, we deduce
$$d_{1+\g}(\sum_{i\in\Z}b_id_{-1+\a_i}v_{a_0-\a_i})
=\sum_{i\in\Z}b_i(-2+\a_i-\g)d_{\g+\a_i}v_{a_0-\a_i}$$
$$=\sum_{i\in\Z}b_i(-2+\a_i-\g)(a_0-\a_i+\a+\b(\g+\a_i))v_{a_0+\g}=0,$$
to yield $\sum_{i\in\Z}b_id_{-1+\a_i}v_{a_0-\a_i}=0$, which is (3.6) for $n=0$.

\medskip
Suppose our claim holds for any $n\le m$ for a fixed $m\in\Z_+$.
Now we consider the claim for $n=m+1$. To get (3.7) for $n=m+1$,
we assume $(3.5')$ for $n=m+1$. For any $\g\in M_0$,
 we compute
$$d_{1+\g}(\sum_{i\in\Z}b_id_{\a_i}d_{-1+a_{m+1}}...d_{-1+a_1}v_{a_0-\a_i})$$
$$=\sum_{i\in\Z}b_i(-1-\g+\a_i)d_{1+\g+\a_i}d_{-1+a_{m+1}}...d_{-1+a_1}v_{a_0-\a_i}
$$
$$+\sum_{i\in\Z}b_id_{\a_i}d_{1+\g}d_{-1+a_{m+1}}...d_{-1+a_1}v_{a_0-\a_i}.
\eqno(3.8)$$ Since
$d_{1+\g}d_{-1+a_{m+1}}...d_{-1+a_1}v_{a_0-\a_i}\in
(L_{-1})^mv_{a_0-\a_i}$  and the coefficients are independent of
$\a_i$, by inductive hypothesis and $(3.5')$, we see that the
second term on the right  hand side of (3.8) is $0$. So
$$d_{1+\g}(\sum_{i\in\Z}b_id_{\a_i}d_{-1+a_{m+1}}...d_{-1+a_1}v_{a_0-\a_i})$$
$$=\sum_{i\in\Z}b_i(-1-\g+\a_i)d_{1+\g+\a_i}d_{-1+a_{m+1}}...d_{-1+a_1}v_{a_0-\a_i}
$$
$$=\sum_{i\in\Z}b_i(-1-\g+\a_i)(-2-\g-\a_i+a_{m+1})d_{\g+\a_i+a_{m+1}}
d_{-1+a_{m}}...d_{-1+a_1}v_{a_0-\a_i}$$
$$
+...+d_{-1+a_{m+1}}...d_{-1+a_{2}}
\sum_{i\in\Z}b_i(-1-\g+\a_i)(-2-\g-\a_i+a_{1})
d_{\g+\a_i+a_{1}}v_{a_0-\a_i}.\eqno(3.9)$$ From inductive
hypothesis and $(3.5')$, we know that each sum on the right hand
side of (3.9) is $0$. Thus, for all $\g\in M_0$,
$$d_{1+\g}(\sum_{i\in\Z}b_id_{\a_i}d_{-1+a_{m+1}}...d_{-1+a_1}v_{a_0-\a_i})=0,$$
which gives
$$\sum_{i\in\Z}b_id_{\a_i}d_{-1+a_{m+1}}...d_{-1+a_1}v_{a_0-\a_i}=0.$$
So (3.7) holds for $n=m+1$.

To verify (3.6) for $n=m+1$, we suppose (3.5) holds for $n=m+1$.
By using (3.7) for $n=m+1$, for any $\g\in M_0$ we deduce that
$$d_{1+\g}\sum_{i\in\Z}b_id_{-1+\a_i}d_{-1+a_{m+1}}...d_{-1+a_1}v_{a_0-\a_i}$$
$$=\sum_{i\in\Z}b_i(-2+\a_i-\g)d_{\g+\a_i}d_{-1+a_{m+1}}
...d_{-1+a_1}v_{a_0-\a_i}
=0,$$
which implies (3.6) for $n=m+1$. Thus Claim 1 holds for $n=m+1$. By inductive principle,
Claim 1 follows.

\vskip 5pt Fix $p\in M_0\setminus\{0\}$. Let $P_i=\{kp|0\le
k\le2i\}$ for any $i\in\N$. For any $n\in \N$, $a\in M_0,
\a_1,...\a_{n}\in \Z p$ with $\a_i\in P_i$, let
$$W(a,\a_1,\a_2,...,\a_{n})\hskip 8cm$$
$$=span
\left\{d_{-1+\a_{n+1}}d_{-1+\a_{n}}...d_{-1+\a_1}v_{a-\a_1-\a_2-...-\a_{n+1}}
 \,|\,
\a_{n+1}\in P_{n+1}\right\}.$$
\medskip
{\bf Claim 2}. {\it For any $\be\in M_0$, we have }
$$d_{-1+\be}d_{-1+a_{n}}...d_{-1+a_1}v_{a-\a_1-...-\a_{n}-\be}
\in W(a,\a_1,\a_2,...,\a_{n}).$$

\medskip
It is clear that we can find nontrivial $b_{\a_{n+1}}\in \bC$ for
$\a_{n+1}\in P_{n+1}$ (consider the following as a linear system
of $2n+3$ unknowns with $2n+3$ equations) such that:
$$\be^k+\sum_{\a_{n+1}\in P_{n+1}}
\a_n^kb_{\a_{n+1}}=0,\f\,\,0\le k\le 2n+2,$$ which is an equality
of (3.5). Then applying (3.6), we see that Claim 2 follows. \vskip
5pt

From Claim 2 we know that
$$
\dim \bar
V_{-n+\a+a}=\Big\{d_{-1+\a_n}d_{-1+a_{n-1}}...
d_{-1+a_1}v_{a-\a_1-\a_2-...-\a_n}:$$
$$
\a_i\in\{kp|0\le k\le 2i\} \Big\}\le  1\cdot3\cdot...\cdot( 2n+1),
\f n\in\Z_+, a\in M_0.
$$
Note that $\dim \bar V_{\a+a}=1$ for all $a\in M_0$.
Thus our theorem follows.\qed

\vskip 5pt It is clear that $\bar V(\a,\b)$ is an irreducible
$\Vir[M]$ module if and only if $ V(\a,\b)$ defined by (3.4) is an
irreducible $\Vir[M_0]$ module. Thus, if we start with the
irreducible $\Vir[M_0]$ module $ V'(\a,\b)$ (which is the
irreducible subquotient of $ V(\a,\b)$ and some of which are no
longer exp-polynomial modules), instead of $ V(\a,\b)$, we get an
irreducible $\Vir[M]$ module $\bar V'(\a,\b)$, which can be also
realized by taking the irreducible subquotient of $\bar V(\a,\b)$
for all $\a,\b\in\bC$.

\vskip 5pt The module $\bar V(\a,\b)$ contains height weight
$\Vir[\Z]$ (the classical Virasoro algebra) modules
$U(\Vir[\Z])v_a$ for any $a\in M_0$. Thus not all its weight
multiplicities are $1$, which indicates the modules $\bar
V(\a,\b)$ are not modules of intermediate series.

%Thus Su's classification theorem is completely wrong.

\vskip 5pt It is natural to ask the following questions: Are
$\Vir[\Z]$ modules \break $U(\Vir[\Z])v_a$  irreducible ? Is it
true that $U(\Vir[\Z])v_a= \oplus_{i\in\Z_+}\bar V_{i+a}$ ? Are
there any other types of irreducible weight modules with finite
dimensional weight spaces over Vir$[M]$ ?

\vskip 5pt It is important to calculate the character formula for
the modules $\bar V(\a,\b)$. We know the following about $\dim
\bar V_{-i+\a+a}$ for $i\in\Z_+, a\in M_0$.

\vskip 5pt {\bf Corollary 3.2.} {\it For any $i\in\Z_+, a_1,
a_2\in M_0$, if $(-i+a_1)(-i+a_2)\ne0$, we have
$$\dim \bar V_{-i+\a+a_1}=\dim \bar V_{-i+\a+a_2}.\eqno(3.10)$$}
{\it Proof.} Suppose $a_1\ne a_2$. Let $a=a_1-a_2$,
$G=\Vir[\Z a]=\oplus_{i\in\Z}\bC d_{ia}$. Then $G$ is isomorphic to the
classical centerless Virasoro algebra. Consider the $G$ module
$W=\oplus_{x\in \Z a}\bar V_{-i+\a+a_2+x}$. From Theorem 3.1, we know that
$W$ is a uniformly bounded $G$ module. By a well-known result, we see that
all the dimensions of $\bar V_{-i+\a+a_2+x}$ are equal except for $\bar V_0$.
The corollary follows.\qed

\bigskip
\vskip.3cm \Refs\nofrills{\bf REFERENCES}
\bigskip
\parindent=0.45in

\leftitem{[BB]} S. Berman, Y. Billig, Irreducible representations for toroidal
Lie algebras, {\it J. Algebra}, 221(1999), no.1, 188--231.

\leftitem{[BS]} S. Berman, J. Szmigielski, Principal realization for the extended
affine Lie algebra of type sl${}\sb 2$ with
coordinates in a simple quantum torus with two generators,
{\it Recent developments in quantum affine algebras and related topics}
(Raleigh, NC, 1998), 39--67, Contemp. Math.,
248, Amer. Math. Soc., Providence, RI, 1999.

\leftitem{[B1]} Y. Billig,
Principal vertex operator representations for toroidal Lie algebras,
 {\it J. Math. Phys.}, 39(1998), No.7, 3844-3864.

\leftitem{[B2]} Y. Billig, Energy-momentum tensor for the toroidal Lie algebras,
preprint
\break
math.RT/0201313.

\leftitem{[C]} V. Chari,  Integrable representations of affine
Lie-algebras, {\it Invent. Math.}, 85(1986), No.2, 317-335.

\leftitem{[DZ]}  D.Z. Djokovic and K. Zhao, Derivations, isomorphisms, and
second cohomology of {generalized} Witt algebras, {\it Trans. of
Amer. Math. Soc.}, Vol.350(2), 643-664(1998).

\leftitem{[E]} S. Eswara Rao, Classification of irreducible integrable modules
for toroidal
Lie algebras with finite dimensional weight spaces,  {\it J. Algebra}, to appear.

\leftitem{[EM]} S. Eswara Rao, R.V. Moody,
 Vertex representations for $n$-toroidal Lie algebras and
 a generalization of the Virasoro algebra, {\it Comm.
  Math. Phys.}, 159(1994), No.2, 239-264.

\leftitem{[G]} Y. Gao, Representations of extended affine Lie algebras
coordinatized by certain quantum tori.
{\it Compositio Math.}, 123 (2000), no. 1, 1-25.

\leftitem{[GL]} M. Golenishcheva-Kutuzova, D.  Lebedev,
Vertex operator representation of some quantum tori Lie algebras, {\it
Comm. Math. Phys.}, 148(1992), no.2, 403-416.

\leftitem{[HWZ]} J. Hu, X. Wang, K. Zhao, Verma modules over generalized
Virasoro algebras Vir$[G]$, {\it J. Pure Appl. Algebra},
177(1), 61-69(2003).

%\leftitem{[KR]} V. G. Kac and Raina, ``Bombay lectures on highest weight
%representations of infinite dimensional Lie algebras,''
%World Sci., Singapore, 1987.

%\leftitem{[OZ]} J.M. Osborn, K. Zhao, Infinite dimensional Lie
%algebras of generalized Block type, {\it Proc. Amer. Math. Soc.},
%Vol.127, No.6, 1641-1650(1999).

\leftitem{[L]} T. A. Larsson, Lowest-Energy Representations of
Non-Centrally Extended Diffeomorphism Algebras, {\it Comm. Math.
Phys.}, {\bf 201}(1999), 461--470.

\leftitem{[M]} O. Mathieu, Classification of Harish-Chandra
modules over the Virasoro algebra, {\it Invent. Math.}, {\bf
107}(1992), 225--234.

\leftitem{[Ma]} V. Mazorchuk, Classification of simple
Harish-Chandra modules over ${Q}$-Virasoro algebra, {\it Math.
Nachr.}, {\bf 209}(2000), 171--177.

\leftitem{[PZ]} J. Patera and H. Zassenhaus, The higher rank
Virasoro algebras, {\it Comm. Math. Phys.}, {\bf 136}(1991),
1--14.

%\leftitem{[R]} R. Ree, On generalized itt algebras,
%{\it Trans. Amer. Math. Soc.} {\bf 83} (1956), 510--546.

%\leftitem{[S]} Y.  Su, Classification
%of Harish-Chandra modules over the higher rank Virasoro
%and super-Virasoro algebras, {\it Comm. Math. Phys.}, (2003).

%\leftitem{[St]} H. Strade, Representations of the Witt algebra,
%{\it J. Alg.} {\bf 49} (2) (1977), 595--605.

\leftitem{[SZ]} Y.  Su and K. Zhao, Generalized Virasoro and
super-Virasoro algebras and modules of intermediate series, {\it
J. Alg.},  252(1), 1-19(2002).

%\leftitem{[W]} R. L. Wilson, Classification of generalized Witt
%algebras over algebraically fields, {\it Trans. Amer. Math. Soc.}
%{\bf 153} (1971), 191--210.

%\leftitem{[Z]} H. Zassenhaus, Ueber Lie'sche Ringe mit
%Primzahlcharakteristik, {\it Hamb. Abh.} {\bf 13} (1939), 1-100.

\endRefs
\vfill
\enddocument
\end